\documentclass [12pt,a4paper,reqno]{amsart}
\input amssymb.sty
\usepackage{amsmath,amsthm,amscd}
\usepackage{amsfonts,amssymb}
\allowdisplaybreaks

\newcommand{\be}{\begin{equation}}
\newcommand{\ef}{\end{equation}}
\chardef\bslash=`\\ 

\newtheorem{thm}{Theorem}[section]
\newtheorem*{thm*}{Theorem}
\newtheorem{cor}[thm]{Corollary}
\newtheorem{lem}[thm]{Lemma}
\newtheorem{prop}[thm]{Proposition}

\theoremstyle{definition}

\newtheorem*{remark*}{Remarks}
\newtheorem*{defn*}{Definition}
\theoremstyle{remark}
\theoremstyle{remark*}

\numberwithin{equation}{section}

\newcommand{\G}{\Gamma}
\newcommand{\wt}{\widetilde}
\newcommand{\wh}{\widehat}
\newcommand{\bk}{\bigskip}
 \renewcommand{\sectionmark}[1]{}

\newcommand{\iy}{\infty}
\newcommand{\Be}{Beltrami}
\newcommand{\hol} {holomorphic}
\newcommand{\qc} {quasiconformal}

\newcommand{\sh} {subharmonic}
\newcommand{\psh} {plurisubharmonic}
\newcommand{\ve}{\varepsilon}

\newcommand{\fc}{\frac}
\newcommand{\Te} {Teichm\"{u}ller}
\newcommand{\const}{\operatorname{const}}

 \usepackage{amsfonts}
\newcommand{\field}[1]{\mathbb{#1}}
\newcommand{\g}{\gamma}

\newcommand{\dl}{\delta}
\newcommand{\D}{\field{D}}
\newcommand{\om}{\omega}
\newcommand{\z}{\zeta}
\newcommand{\ov}{\overline}
\newcommand{\vp}{\varphi}
\newcommand{\hC}{\widehat{\field{C}}}
\newcommand{\C}{\field{C}}

\newcommand{\B}{\mathbf{B}}
\newcommand{\T}{\mathbf{T}}

\newcommand{\Hol}{\operatorname{Hol}}
\newcommand{\id}{\operatorname{id}}
\newcommand{\grad}{\operatorname{grad}}

\newcommand{\uTs} {universal Teichm\"{u}ller space}
\newcommand{\Belt} {\operatorname{Belt}}
\newcommand{\Ca} {Carath\'{e}odory}
\newcommand{\Gr} {Grunsky}
\newcommand{\Ko} {Kobayashi}
\newcommand{\Om} {\Omega}
\newcommand{\vk} {\varkappa}
\newcommand{\kp} {\kappa}
\newcommand{\x} {\mathbf x}
\renewcommand{\a} {\alpha}
\newcommand{\ld}{\lambda}
\newcommand{\Ld}{\Lambda}

\newcommand{\dist}{\operatorname{dist}}

\begin{document}

\title{Milin's coefficients, complex geometry of Teichm\"{u}ller spaces
and variational calculus for univalent functions}

\author{Samuel L. Krushkal}

\begin{abstract} We investigate the invariant metrics and complex
geodesics in the \uTs \ and \Te \ space of the punctured disk using
Milin's coefficient inequalities. This technique allows us to
establish that all non-expanding invariant metrics in either of
these spaces coincide with its intrinsic Teichm\"{u}ller metric.

Other applications concern the variational theory for univalent
functions with quasiconformal extension. It turns out that geometric
features caused by the equality of metrics and connection with
complex geodesics provide deep distortion results for various
classes of such functions and create new phenomena which do not
appear in the classical geometric function theory.
\end{abstract}

\date{\today\hskip4mm ({\tt geovar.tex})}

\maketitle

\bigskip

{\small {\textbf {2010 Mathematics Subject Classification:} Primary:
30C55, 30C62, 30C75, 30F60, 32F45; Secondary: 30F45, 46G20}

\medskip

\bigskip

\markboth{Samuel Krushkal}{Complex geometry and variational
calculus} \pagestyle{headings}

\bk
\section{Key theorems on invariant metrics and geodesics}

\subsection{Preamble}
The Milin coefficient inequalities arose as a generalization of the
classical Grunsky inequalities but coincide with the later only for
conformal maps of the unit disk.

We apply a quasiconformal variant of these inequalities to
investigation of complex metric geometry and complex geodesics on
two \Te \ spaces:  the universal space and \Te \ space of the
punctured disk and apply their geometry to variational calculus for
univalent functions on the generic quasidisks with \qc \ extensions.
Such functions play an important role in the theory of \Te \ spaces
and also form one of the basic classes in geometric function theory.

It will be shown that the intrinsic geometric features provide deep
distortion results, in particular, allow one to solve explicitly
some general variational problems. On the other hand, they cause
surprising phenomena which do not arise in the classical variational
theory for univalent functions.

\subsection{Main property of invariant metrics of \Te \ spaces}
We shall use the notations $\hC = \C \cup \{\iy\}, \ \D = \{|z| <
1\}, \ \D^* = \hC \setminus \ov{\D} = \{|z| > 1\}$ and consider two
\Te \ spaces: the \uTs \ $\T = \T(\D)$ and the \Te \ space $\T_1 =
\T(\D^0)$ of the punctured disk $\D^0 = \D \setminus \{0\}$ endowed
with the homotopy class of \qc \ homeomorphisms containing the
identity map and regarded as the base point of $\T(\D^0)$. The space
$\T_1$ is model for \Te \ spaces of punctured disks with arbitrary
number of punctures and even for more general flat Riemann surfaces.

Every \Te \ space $\wt \T$ is a complex Banach manifold, thus it
possesses the invariant Carath\'{e}\-odory and \Ko \ distances (the
smallest and the largest among all holomorphically non-expanding
metrics). Denote these metrics by $c_{\wt \T}$ and $d_{\wt \T}$, and
let $\tau_{\wt \T}$ be the intrinsic \Te \ metric of this space
canonically determined by \qc \ maps. The corresponding
infinitesimal Finsler metrics (defined on the tangent bundle
$\mathcal T \wt \T$ of $\wt \T$) are denoted by $\mathcal C_{\wt
\T}$ and $\mathcal K_{\wt \T}$ and $F_{\wt \T}$, respectively. Then
 \be\label{1.1}
c_{\wt \T}(\cdot, \cdot) \le d_{\wt \T}(\cdot, \cdot) \le \tau_{\wt
\T}(\cdot, \cdot),
\end{equation}
and by the Royden-Gardiner theorem the metrics $d_{\wt \T}$ and
$\tau_{\wt \T}$ (and their infinitesimal forms) are equal, see, e.g.
\cite{EKK}, \cite{EM}, \cite{GL}, \cite{Ro}.

In view of applications, we mainly focus on the \Ca \ metric of the
space $\T_1$ and first establish that it equals the \Te \ metric.
This yields the all non-expanding invariant metrics on $\T_1$ agree
with $\tau_\T$, and the \Te \ extremal disks are geodesic with
respect to all invariant metrics.

\begin{thm} The \Ca \ metric of the space $\T_1$ coincides with
its \Ko \ metric, hence all invariant non-expanding metrics on
$\T_1$ are equal its \Te \ metric, and
 \be\label{1.2}
 c_{\T_1}(\vp, \psi) = d_{\T_1}(\vp,\psi) = \tau_{\T_1}(\vp, \psi) = \inf \{
d_\D(h^{-1}(\vp), h^{-1}(\psi)): \ h \in \Hol(\D, \T_1)\},
\end{equation}
where $d_\D$ denotes the hyperbolic metric of the unit disk of
curvature $-4$.

Similarly, the infinitesimal forms of these metrics coincide with
the Finsler metric $F_{\T_1}(\vp, v)$ generating $\tau_{\T_1}$ and
have \hol \ sectional curvature $-4$.
\end{thm}

Such a result is known only for the \uTs \ and underlies various
applications; its proof was given in \cite{Kr4} (and somewhat
modified in \cite{Kr7}). In view of importance, we present this fact
here as a separate theorem giving its simplified proof and new
applications.

\begin{thm} All invariant non-expanding metrics on the \uTs \
$\T$ are equal to its \Te \ metric.
\end{thm}

The proof of both theorems involves the Grunsky-Milin coefficient
inequalities.

\subsection{Complex geodesics}
The equality of metrics allows one to describe complex geodesics in
the spaces $\T$ and $\T_1$. Let $\wt \T$ denote either of these
spaces.

Recall that if $X$ is a domain in a complex Banach space $E$ endowed
with a pseudo-distance $\rho_X$, then a \hol \ map $h: \ \D \to X$
is called a {\bf complex $\rho$-geodesic} if there exist $t_1 \ne
t_2$ in $\D$ such that
$$
d_\D(t_1, t_2) = \rho_X(h(t_1), h(t_2));
$$
one says also that the points $h(t_1)$ and $h(t_2)$ can be
joined by a complex $\rho$-geodesic (see \cite{Ve})).

If $h$ is a complex $c_X$-geodesic then it also is $d_X$-geodesic
and the above equality holds for all points $t_1, \ t_2 \in \D$, so
$h(\D)$ is a \hol \ disk in $X$ hyperbolically isometric to $\D$. As
an important consequence of Theorems 1.1. and 1.2, one gets the
following result where the complex geodesics are understanding in
the strongest sense, i.e., as $c_{\wt T}$-geodesics.

\bk
\begin{thm} (i) Any two points of the space $\wt \T$ can be joined by
a complex geodesic. The geodesic joining a Strebel's point with the
base point is unique and defines the corresponding \Te \ extremal
disk.

(ii) For any point $\vp \in \wt \T$ and any nonzero tangent vector
$v$ at this point, there exists at least one complex geodesic $h: \
\D \to \wt \T$ such that $h(0) = \vp$ and $h^\prime(0)$ is collinear
to $v$.
\end{thm}

\subsection{Geometric and analytic features}
The following consequence of Theorem 1.1 relates to pluripotential
features of $\wt \T$ and is useful in variational problems on
compact subsets of $\Sigma^0(D)$.

\begin{cor} Any non-expanding invariant metrics $\rho$ on
the space $\T_1$ with the base point $D^0$ relates to the similar
metric $\rho_{B_k}$ on hyperbolic balls $B_k = \{\psi \in \T_1: \
\tau_\T(\psi, \mathbf 0) < \tanh k\} \ \ (0 < k < 1)$ by
$$
\rho_{B_k}(\psi_1, \psi_2) = \tanh^{-1}
\Bigl(\frac{l(\rho_\T(\psi_1, \psi_2))}{k}\Bigr) = d_\D \Bigl(0,
\frac{l(d_\T(\psi_1, \psi_2))}{k}\Bigr), \quad l(s) = \tanh s.
$$
Similar relation holds for the pluricomplex Green functions of the
space $\T_1$ and its balls $B_k(\T_1)$.

For the \uTs \ $\T$, this was established in  \cite{Kr5}. The proof
for the space $\T_1$ follows the same lines using Theorem 1.1.
\end{cor}

This assertion is obtained from Theorem 1.1 using the arguments
applied in \cite{Kr5} for the \Ko \ metric of \uTs.

It is not known, how to relate the invariant distances of the balls
in generic complex manifolds $X$ with the corresponding distances on
$X$.

\bk The following corollary controls the growth of \hol \ maps of
$\wt \T$ on geodesic disks.

\begin{cor} If a \hol \ map  $J: \ \wt \T \to \D$ into the unit disk
is such that its restriction to a geodesic disk $\D(\mu_0) =
\{\phi_{\wt \T}(t \mu_0/\|\mu_0\|_\iy ): \ |t| < 1\}$ has at the
origin zero of order $m$, i.e.,
 \be\label{1.3}
J_{\mu_0}(t) := J \circ \phi_{\wt \T}(t \mu_0/\|\mu_0\|_\iy) = c_m
t^m + c_{m+1} t^{m+1} + \dots,
\end{equation}
then the growth of $|J|$ on this disk is estimated by
 \be\label{1.4}
|J_{\mu_0}(t)| \le \tanh \Bigl(|t|^m \fc{|t| + |c_m|}{1 +
|c_m||t|}\Bigr) \le d_{\wt \T} \Bigl(\mathbf 0, \phi_{\wt \T}
\Bigl(t^m \fc{\mu_0}{\|\mu_0\|_\iy}\Bigr)\Bigr).
\end{equation}
The equality in the right inequality occurs (even for one $t_0 \ne
0$) only when $|c_m| = 1$; then $J_{\mu_0}(t)$ is a hyperbolic
isometry of the unit disk and all terms in (1.3) are equal.
\end{cor}

\subsection{} The above theorems and corollaries have deep applications
to geometric function theory. Some of those are presented in the
last two sections.

\section{Background}
We recall some notions and results which will be used in the proofs
of the above theorems.

\subsection{Invariant metrics on \Te \ spaces}
Let $L$ be a bounded oriented quasicircle in the complex plane $\C$
with the interior and exterior domains $D$ and $D^*$ so that $D^*$
contains the infinite point $z = \iy$. Consider the unit ball of \Be
\ coefficients supported on $D$,
$$
\Belt(D)_1 = \{\mu \in L_\iy(C): \ \mu(z)|D^* = 0, \ \ \|\mu\|_\iy <
1\}
$$
and their pairing with $\psi \in L_1(D)$ by
$$
\langle \mu, \psi\rangle_D = \iint\limits_D  \mu(z) \psi(z) dx dy
\quad (z = x + i y).
$$
The following two sets of \hol \ functions $\psi$ (equivalently, of
\hol \ quadratic differentials $\psi dz^2$)
$$
\begin{aligned}
A_1(D) &= \{\psi \in L_1(D): \ \psi \ \ \text{\hol \ in} \ \ D\},  \\
A_1^2(D) &= \{\psi = \om^2 \in A_1(D): \ \om \ \ \text{\hol \ in} \
\ D\}
\end{aligned}
$$
are intrinsically connected with the extremal \Be \ coefficients
(hence, with the \Te \ norm) and Grunsky-Milin inequalities.

The elements of $A_1^2$ can be regarded as the squares of abelian
\hol \ differentials on $D$.

Rescaling the domain $D^*$ to have $D^* = f^{\mu_0}(\D^*)$ for some
$f^{\mu_0}(z) = z + b_0 + b_1 z^{-1} + \dots$ preserving $z = 0$
with $\mu_0 \in \Belt(\D)_1$, one can use this domain as a new base
point of the \uTs \ $\T$ whose points are the equivalence classes
$[\mu]$ of $\mu \in \Belt(D)_1$ so that
$$
\mu_1 \sim \mu_2 \ \ \text{if} \ \ w^{\mu_1}(z) = w^{\mu_2}(z) =
f(z) \ \ \text{on} \ \ \ov{D^*}.
$$
We shall also denote such classes by $[f]$. This space is modeled as
a bounded domain in the complex Banach space $\B(D^*)$ of the
\textbf{Schwarzian derivatives}
$$
S_w = (w^{\prime\prime}/w^\prime)^\prime -
(w^{\prime\prime}/w^\prime)^2/2, \quad w = f^\mu|D^*,
$$
of locally univalent functions on $D_0^*$ with norm $ \|\vp\| =
\sup_{D^*} \ld_{D^*}(z)^{-2} |\vp(z)|$, where $\ld_{D^*}(z) |dz|$ is
the differential hyperbolic metric on $D^*$ of curvature $-4$. This
modeling domain is filled by the Schwarzians of globally univalent
functions on $D^*$ with \qc \ extension.

For the unit disk $\D$, $\ld_\D(z) = 1/(1 - |z|^2)$, and the global
hyperbolic distance
$$
d_\D(z_1, z_2) = \tanh^{-1} [(z_1 - z_2)/(1 - \ov z_1 z_2)].
$$

The intrinsic \textbf{\Te \ metric} of the space $\T$ is defined by
$$
\tau_\T (\phi_\T (\mu), \phi_\T (\nu)) = \frac{1}{2} \inf \bigl\{
\log K \bigl( w^{\mu_*} \circ \bigl(w^{\nu_*} \bigr)^{-1} \bigr) : \
\mu_* \in \phi_\T(\mu), \nu_* \in \phi_\T(\nu) \bigr\},
$$
where $\phi_\T$ is the factorizing \hol \ projection $\Belt(D)_1 \to
\T$. This metric is the integral form of the infinitesimal Finsler
metric (structure)
 \be\label{2.1}
  F_\T(\phi_\T(\mu), \phi_\T^\prime(\mu) \nu) = \inf
\{\|\nu_*/(1 - |\mu|^2)^{-1}\|_\iy: \ \phi_\T^\prime(\mu) \nu_* =
\phi_\T^\prime(\mu) \nu\}
\end{equation}
on the tangent bundle $\mathcal T\T$ of $\T$, which is locally
Lipschitzian (see \cite{EE}).

Note also that $\tau_\T(\mathbf 0, S_f) = \tanh^{-1} k(f)$, where
$k(f)$ is the \textbf{\Te \ norm} of a univalent function $f$.

The \textbf{\Ko} \ and \textbf{\Ca} \ metrics $d_{\wt \T}$ and
$c_{\wt \T}$ of a \Te \ space $\wt \T$ relate to complex structure
of this space and are defined, respectively, as the largest
pseudometric $d$ on $\wt \T$ which does not get increased by the
\hol \ maps $h: \ \D \to \wt \T$ so that for any two points $\vp_1,
\ \vp_2 \in \wt \T$, we have
$$
d_{\wt \T}(\vp_1, \vp_2) \leq \inf \{d_\D(0,t): \ h(0) = \vp_1, \
h(t) = \vp_2\},
$$
and
$$
c_{\wt \T}(\vp_1, \vp_2) = \sup d_\D(h(\vp_1), h(\vp_2)),
$$
taking the supremum over all  \hol \ maps $h: \ \wt \T \to \D$.

The corresponding infinitesimal forms of the \Ko \ and \Ca \ metrics
are defined for the points $(\vp, v) \in \mathcal T \wt \T$,
respectively, by
$$
\begin{aligned}
\mathcal K_{\wt \T}(\vp, v) &= \inf \{1/r: \ r > 0, \ h \in
\Hol(\D_r, \wt \T),
h(0) = \vp, h^\prime (0) = v\},  \\
\mathcal C_{\wt \T}(\vp, v) &= \sup \{|d f(\vp) v|: \ f \in \Hol(\wt
\T, \D), f(\vp) = 0\},
\end{aligned}
$$
where $\Hol(X, Y)$ denotes the collection of \hol \ maps of a
complex manifold $X$ into $Y$ and $\D_r$ is the disk $\{|z| < r\}$.

The sectional \textbf{\hol \ curvature} $\kp_F(x, v)$ of a Finsler
metric $F(x, v)$ on (the tangent bundle of) a complex Banach
manifold $X$ is defined as the supremum of the generalized Gaussian
curvatures
 $$
 \kappa_\ld(t) = - \fc{\Delta \log \ld(t)}{\ld(t)^2}
 \quad \text{for} \ \ \ld(t) = F(h(t), h^\prime(t))
 $$
over appropriate collections of \hol \ maps $h$ from the disk into
$X$ for a given tangent direction $v$ in the image. Here $\Delta$
means the generalized Laplacian
$$
\Delta \ld(t) = 4 \liminf\limits_{r \to 0} \frac{1}{r^2} \Big\{
\frac{1}{2 \pi} \int_0^{2\pi} \ld(t + re^{i \theta}) d \theta -
\ld(t) \Big\}
$$
(provided that $0 \le \ld(t) < \iy$). Similar to $C^2$ functions,
for which $\Delta$ coincides with the usual Laplacian $4
\ov{\partial}\partial$, one obtains that $\ld$ is \sh \ on a domain
$\Om$ if and only if $\Delta \ld(t) \ge 0$; hence, at the points
$t_0$ of local maxima of $\ld$ with $\ld(t_0) > - \iy$, we have
$\Delta \ld(t_0) \le 0$.

Generically, the \hol \ curvature of the \Ko \ metric $\mathcal
K_X(x, v)$ of any complete hyperbolic manifold $X$ satisfies
$\kp_{\mathcal K_X}(x, v) \ge - 4$ at all points $(x, v)$ of the
tangent bundle $\mathcal T(X)$ of $X$, and for the \Ca \ metric
$\mathcal C_X$ we have $\kp_{\mathcal C_X}(x, v) \le - 4$.

\subsection{The Grunsky and Milin coefficients inequalities}
Denote by $\Sigma(D^*)$ the collection of univalent functions $f$ in
a quasidisk $D^*$ with hydrodynamical expansion
 \be\label{2.2}
f(z) = z + b_0 + b_1 z^{-1} + \dots \quad \text{near} \ \ z = \iy,
 \end{equation}
and let $\Sigma^0(D^*)$ denote its subset formed by functions having
\qc \ extensions across the boundary (hence to $\C$). Each $f \in
\Sigma(D^*)$ determines a \hol \ map
 \be\label{2.3}
 - \log \fc{f(z) - f(\z)}{z - \z}
= \sum\limits_{m, n = 1}^\iy \fc{\a_{m n}}{\chi(z)^m \ \chi(\z)^n}:
\ (D_0^*)^2 \to \hC
\end{equation}
where $\chi$ is the conformal map of $D^*$ onto the disk $\D^*$ with
$\chi(\iy) = \iy, \ \chi^\prime(\iy) > 0$, and the Taylor
coefficients $\a_{m n}$ are called the Milin coefficients of $f$. In
the classical case $D^* = \D^*$, those are the standard \Gr \
coefficients.

Due to the \Gr \ univalence theorem \cite{Gr} and its Milin's
extension \cite{Mi}, a function $f$ \hol \ near the infinity (with
hydrodynamical normalization) is extended to a univalent function in
the whole domain $D^*$ if and only if its coefficients $\a_{m n}$
satisfy the inequality
$$
\Big\vert \sum\limits_{m,n = 1}^{\iy} \ \sqrt{m n} \ \a_{mn} \ x_m
x_n \Big\vert \le 1
$$
for any sequence $\x = (x_n) \in l^2$ with $\|\x \|^2 =
\sum\limits_1^\iy |x_n|^2 = 1$. We denote the unit sphere of these
Hilbert space by $S(l^2)$ and call the quantity
  \be\label{2.4}
\vk_{D^*}(f) := \sup \Big\{ \Big\vert \sum\limits_{m,n = 1}^{\iy} \
\sqrt{m n} \ \a_{mn} \ x_m x_n \Big\vert : \ {\mathbf x} = (x_n) \in
S(l^2)\Big\}
\end{equation}
the \textbf{Grunsky norm} of $f$ on $D^*$. For $D^* = \D^*$, we
shall use simplified notations $\Sigma$ and $\vk(f)$.

Noting that each coefficient $\a_{m n}(f)$ in (1.2) is represented
as a polynomial of a finite number of the initial coefficients $b_1,
b_2, \dots, b_{m+n-1}$ of $f$, one derives after normalizing \qc \
extensions of $f^\mu$ in $D$ (for example, by $f(0) = 0$) the \hol \
dependence of $\beta_{m n}(f)$ on \Be \ coefficients  $\mu$ and on
the Schwarzian derivatives $S_f$ on $D^*$ runing over the \uTs \
$\T$ with the base point $\chi^\prime(\iy) D^*$.

For any finite $M, N$ and $1 \le j \le M, \ 1 \le l \le N$, we have
$$
\Big\vert \sum\limits_{m=j}^M \sum\limits_{n=l}^N \ \sqrt{m n} \
\a_{mn} x_m x_n \Big\vert^2 \le \sum\limits_{m=j}^M |x_m|^2
\sum\limits_{n=l}^N |x_n|^2;
$$
this inequality is a consequence of Milin's univalence theorems (cf.
[Mi, p. 193], [Po, p. 61]). Thus for each $\x = (x_n) \in S(l^2)$,
the function
 \be\label{2.5} h_{\x}(\mu) =
\sum\limits_{m,n=1}^\iy \ \sqrt{m n} \ \a_{m n} (f^\mu) x_m x_n
\end{equation}
maps the space $\T$  holomorphically into the unit disk, and
 \be\label{2.6}
\sup_{\x} |h_{\x}(f^\mu)| = \vk_{D^*}(f^\mu).
\end{equation}
This implies also that the \Gr \ norm $\vk_{D^*}(f^\mu)$ is a
continuous \psh \ function on $\Belt(D)_1$ and on the space $\T$
(cf. \cite{Kr2}, \cite{Kr8}).

\bk The following key results obtained in \cite{Kr2}, \cite{Kr8} by
applying the maps (2.6) underly the proofs of Theorems 1.1 and 1.2.

\begin{prop} (a) \ The \Gr \ norm $\vk_{D^*}(f)$ of every function
$f \in \Sigma^0(D^*)$ is estimated by its \Te \ norm $k = k(f)$ by
 \be\label{2.7}
\vk_{D^*}(f) \le k \fc{k + \a_D(f)}{1 + \a_D(f) k},
\end{equation}
where
$$
\a_D(f^\mu) = \sup \ \{|\langle \mu, \vp\rangle_D|: \ \vp \in
A_1^2(D), \ \|\vp\|_{A_1} = 1 \Big\} \le 1.
$$
and $\vk_{D^*}(f) < k$ unless $\a_D(f) = 1$. The last equality
occurs if and only if $\vk_{D^*}(f) = k(f)$.

(b)\ The equality $\vk_{D^*}f = k(f)$ holds if and only if the
function $f$ is the restriction to $D^*$ of a \qc \ self-map
$w^{\mu_0}$ of $\hC$ with \Be \ coefficient $\mu_0$ satisfying the
condition
 \be\label{2.8} \sup |\langle \mu_0, \vp\rangle_D| =
\|\mu_0\|_\iy,
\end{equation}
where the supremum is taken over \hol \ functions $\vp \in A_1^2(D)$
with $\|\vp\|_{A_1} = 1$.

If, in addition, the equivalence class of $f$ (the collection of
maps equal $f$ on $\partial D^*$) is a Strebel point, then $\mu_0$
is necessarily of the form
 \be\label{2.9}
\mu_0(z) = \|\mu_0\|_\iy |\psi_0(z)|/\psi_0(z) \quad \text{with} \ \
\psi_0 \in A_1^2(D).
\end{equation}
\end{prop}

The condition (2.8) has a geometric nature. Its proof in \cite{Kr2},
\cite{Kr8} shows that the functions (2.5) generate a maximizing
sequence on which the \Ca \ distance $c_\T(\mathbf 0,
S_{f^{t\mu_0}})$ is attained and equals the \Te \ distance (compare
with Kra's theorem in \cite{K} about the \Ca \ metric on \Te \
abelian disks for Riemann surfaces with finitely generated
fundamental groups).

In a special case, when the domain $D^*$ is the disk $\D^*$ and $f$
is analytic up to its boundary $\{|z| = 1\}$, the equality (2.9) was
obtained by a different method in \cite{Ku3}.

\bk \noindent{\bf Remark}. The \Gr \ coefficients were originally
defined in \cite{Gr} for finitely connected plane domains and have
been later generalized to bordered Riemann surfaces (see \cite{SS},
\cite{Le}). Milin's coefficients were introduced in \cite{Mi} also
for finitely connected domains. Both types of coefficients coincide
only for conformal maps of a circular disk.

\subsection{Weak$^*$ compactness of \hol \ families in Banach
domains} One of the main underlying facts in the proof of main
theorems is the existence of \hol \ maps from the space $\T$ or
$\T_1$ onto the unit disk on which the \Ca \ distance is attained.
It relies on classical results on compactness in the dual week$^*$
topology. The following two propositions were related to me by David
Shoikhet.

An analog of Montel's theorem for the infinite dimensional case is
given by

\begin{prop} Let $X$ and $Y$ be Banach spaces
and let $G$ be a domain in $X$. A bounded set $\Om$ in $\Hol(G, Y)$
is relatively compact with respect to the topology of uniform
convergence on compact subsets of $D$ (compact open topology on
$\Hol(G, Y)$) if and only if for each $x$ in $G$ the orbit $h(x)$,
where $h$ runs over $\Om$, is relatively compact in $Y$.
\end{prop}

In fact, this proposition is a consequence of the classical Ascoli
theorem. As a consequence of the Alaoglu-Bourbaki theorem, one
derives the following result.

\begin{prop} Let $Y$ be reflexive and let $\Om$ be a bounded set
in  $\Hol(G, Y)$. Then any sequence from $\Om$ contains a
subsequence which weakly converges to a \hol \ map from $G$ to $Y$.
\end{prop}

Hence, the compactness (strong or weak) is actually required only
for the image of $G$ in $Y$ under a given family of maps.

In our case $Y = \C$, and for each $x \in \wt \T$ its orbit $h(x)$
is located in the unit disk which is compact. This implies that for
any point $x_0$ from $\wt \T$ there exists a \hol \ map $h_0: \ \wt
\T \to \D$ with $h_0(0) = 0$ and $d_\D(0, h_0(x_0)) = c_\T(\mathbf
0, x_0)$.

\section{Proof of Theorem 1.2}

It suffices to establish the equality of the \Ca \ and \Te \ metrics
of $\T$ for the equivalence classes $[\mu]$ which are the
\textbf{Strebel points} of $\T$. This means that $[\mu]$ contains
(unique) \Be \ coefficient of the form $\mu_0(z) = k
|\psi_0(z)|/\psi_0(z)$, where $k < 1$ and $\psi_0 \in A_1(\D)$.
Accordingly, we have in $\T$ the \Te \ geodesic disks
$$
\D(\mu_0) = \{\phi_\T(t \mu_0): t \in \D\}.
$$
It is well known that the set of Strebel points is open and dense in
any \Te \ space (see \cite{GL}, \cite{St}); in addition, both
metrics are continuous on this space.

If the defining quadratic differential $\psi_0$ has in $\D$ only
zeros of even order, i.e. belongs to $A_1^2(\D)$, then the
equalities (1.1) follow from Proposition 2.1 with $D_0 = \D$,
because the \hol \ maps (2.5) provide in view of the equality (2.8)
a maximizing sequence for the \Ca \ distance $c_\T(\mathbf 0,
S_{f^{\mu_0}})$ and this distance equals $\tanh^{-1} k(f^{\mu_0})$.

Namely, using the variation of  $f^\mu \in \Sigma^0$ with $f^\mu(0)
= 0$ by
 \be\label{3.1}
 f^\mu(z) = z - \fc{1}{\pi} \iint\limits_\D
\mu(w) \left( \fc{1}{w - z} - \fc{1}{w} \right) du dv +
O(\|\mu\|_\iy^2), \quad w = u + iv
\end{equation}
(with uniformly bounded ratio $O(\|\mu\|_\iy^2)/\|\mu\|_\iy^2$ on
compact subsets of $\C$), one derives that the \Gr \ coefficients of
$f^\mu$ are varied by
 \be\label{3.2}
\a_{mn}(S_{f^\mu}) = - \fc{1}{\pi} \iint\limits_\D \mu(z) z^{m+n-2}
dx dy + O(\|\mu\|_\iy^2), \quad \|\mu\|_\iy \to 0.
\end{equation}
This implies that the differential at zero of the corresponding map
$$
\wh h_\x = h_\x \circ \phi_\T: \Belt(\D)_1 \to \T \to \D \quad
\text{with} \ \ \x = (x_n) \in S(l^2)
$$
is given by
 \be\label{3.3}
 d \wh h_{\x}[\mathbf 0] (\mu/\|\mu\|_\iy) = - \fc{1}{\pi}
\iint\limits_\D \fc{\mu(z)}{\|\mu\|_\iy} \ \sum\limits_{m+n=2}^\iy
\sqrt{m n} \ x_m x_n z^{m+n-2} dx dy.
\end{equation}
On the other hand, as was established in \cite{Kr2}, the elements of
$A_1^2(\D)$ are represented in the form
  \be\label{3.4}
\psi(z) = \om(z)^2 = \fc{1}{\pi} \sum\limits_{m+n=2}^\iy \ \sqrt{m
n} \ x_m x_n z^{m+n-2},
\end{equation}
with $\|\x\|_{l^2} = \|\om\|_{L_2}$. The relations (2.6), (3.3),
(3.4) together with Schwarz's lemma, imply for $\mu = \mu_0$ the
desired equality
  \be\label{3.5}
\tanh c_\T(\mathbf 0, S_{f^{\mu_0}}) = \vk(f^{\mu_0}) = k.
 \end{equation}

The investigation of the generic case, when $\psi_0$ has in $\D$ a
finite or infinite number of zeros of odd order, involves the Milin
coefficient inequalities.

First recall the chain rule for \Be \ coefficients: for any $\mu, \
\nu \in \Belt(\C)_1$, the solutions $w^\mu$ of the corresponding \Be
\ equation $\partial_{\ov z} w = \mu \partial_z w$ on $\hC$ satisfy
$w^\mu \circ w^\nu = w^{\sigma_\nu(\mu)}$, with
 \be\label{3.6}
\sigma_\nu(\mu) = (\nu + \mu^*)/(1 + \ov \nu \mu^*),
\end{equation}
where
$$
\mu^*(z) = \mu \circ w^\nu(z) \ \ov{\partial_z w^\nu(z)}/\partial_z
w^\nu(z).
$$
Thus, for $\nu$ fixed, $\sigma_\nu(\mu)$ depends holomorphically on
$\mu$ as a map $L_\iy(\C) \to L_\iy(\C)$.

Without loss of generality, one can assume that the function
$\psi_0$ does not have at $z = 0$ (hence at some disk $\{|z| < d <
1\}$) zero of odd order. Otherwise, after squaring
$$
f^{\mu_0} \mapsto \mathcal R_2 f^{\mu_0} = f^{\mu_0}(z^2)^{1/2} = z
+ \fc{b_0}{2} \fc{1}{z} + \fc {b_3^\prime}{z^3} + \dots
$$
one obtains an odd function from $\Sigma^0$ whose \Be \ coefficient
on $\D$ equals $\mathcal R_2^* \mu_0 = \mu_0(z^2) \ov z/z$ being
defined by quadratic differential $\mathcal R_2^* \psi_0 = 4
\psi_0(z^2) z^2$ with zero of even order at the origin.

In addition, the Taylor and \Gr coefficients of $\mathcal R_2 f^\mu$
are represented as polynomials of the initial Taylor coefficients
$b_1, \dots \ , b_s$ of the original function $f^\mu$, thus $\a_{m
n}(\mathcal R_2 f^\mu)$, together with $\a_{m n}(f^\mu)$, depend
holomorphically on $\mu$ and $S_{f^\mu}$.

Now fix a $\dl > 0$ close to $1$ and delete from the disk $\D$ the
annulus $\mathcal A_\dl = \{\dl < |z| < 1\}$ and the circular
triangles $\Delta_1, \dots \ , \Delta_{m(\dl)}$ such that the base
of each $\Delta_j$ is an arc of the circle $\{|z| = \dl\}$, its
opposite vertex is a zero $a_j$ of odd order satisfying $\rho \le
|a_j| < \dl$, and two other sides of $\Delta_j$ are the straight
line segments symmetric with respect to the radial segment
$\bigl[a_j, \dl e^{i \arg a_j} \bigr]$. In the case when several
zeros are located on the same radius, it suffices to take only the
zero with minimal modulus. Denote
$$
E_\dl  = \cup_j \Delta_j \cup \mathcal A_\dl, \ \ D_\dl = \D
\setminus \ov{E_\dl} \Subset \D, \ \ D_\dl^* = \D^* \cup E_\dl = \hC
\setminus \ov{D_\dl}
$$
and put
$$
\mu_1(z) = \begin{cases} \mu_0(z)&      \text{if} \ \ z \in E_\dl,   \\
                          0&  \text{otherwise}.
\end{cases}
$$
Note that $f^{\mu_1}$ is normalized by (2.2) and $f^{\mu_1}(0) = 0$,
so $S_{f^{\mu_1}} \in \T$. Letting $\mu_2 = \mu_0 - \mu_1$, one
factorizes the initial automorphism $f^{\mu_0}$ of $\hC$ via
$$
f^{\mu_0} = f^{\sigma_0} \circ f^{\mu_1}
$$
with
 \be\label{3.7}
\sigma_0  = (f^{\mu_1})^* \mu_0 = \Bigl(\fc{\mu_2}{1 - \ov{\mu_1}
\mu_0} \ \fc{\partial_z f^{\mu_1}}{\ov{\partial_z f^{\mu_1}}}\Bigr)
\circ (f^{\mu_1})^{-1} \in \Belt(f^{\mu_1}(D_\dl))_1.
\end{equation}
Since  $f^{\mu_1}$ is conformal on $D_\dl$, the coefficient
$\sigma_0$ is represented by $\sigma_0 = k |\psi_\dl|/\psi_\dl$ with
  \be\label{3.8}
 \psi_\dl(w) = (\psi_0 \circ
\wt f) \ (\wt f^\prime)^2(w) \in A_1^2( f^{\mu_1}(D_\dl)), \quad \wt
f = (f^{\mu_1})^{-1};
\end{equation}
this coefficient is extremal in its class in the ball
$\Belt(f^{\mu_1}(D_\dl))_1$.

The equivalence classes of \Be \ coefficients $\nu \in
\Belt(f^{\mu_1}(D_\dl))_1$ under the relation $\nu_1 \sim \nu_2$ if
$w^{\nu_1} = w^{\nu_2}$ on $\partial f^{\mu_1}(D_\dl)$ form the
quotient space $\T^* = \T(f^{\mu_1}(D_\dl))$ which is
biholomorphically isomorphic to the \uTs \ with the base point
$f^{\mu_1}(D_\dl)$. The factorizing projection $\phi_{\T^*}: \
\Belt(f^{\mu_1}(D_\dl))_1 \to \T^*$ is a \hol \ split submersion,
which means that it has local \hol \ sections.

The chain rule for the Schwarzians
$$
S_{f_2\circ f_1} = (S_{f_2} \circ f_1) (f_1^\prime)^2 + S_{f_1}
$$
applied to $w^\nu \circ f^{\mu_1}$, where $w^\nu \in
\Sigma(D_\dl^*)$ creates a \hol \ map $\eta: \ \T \to \T^*$ moving
the base point to the base point.

Now, applying to $w^\nu \in \Sigma^0(f^{\mu_1}(D_\dl))$ the
variation of type (3.1), one obtains the following generalizations
of (3.2) to Milin's coefficients given in \cite{Kr8}
$$
\a_{m n}(S_{w^\nu}) = - \fc{1}{\pi} \iint\limits_{f^{\mu_1}(D_\dl)}
\ \nu(w) P_m^\prime(w) P_n^\prime(w) du dv + O(\|\nu\|_\iy^2),
$$
and accordingly, instead of (3.3),
$$
d \wh h_\x[\mathbf 0] (\nu/\|\nu\|_\iy) = - \fc{1}{\pi}
\iint\limits_{f^{\mu_1}(\D_\dl)} \ \fc{\nu(z)}{\|\nu\|_\iy} \
\sum\limits_{m,n=1}^\iy
 x_m x_n \ P_m^\prime(z) P_n^\prime(z) dx dy.
 $$
Here $\wh h_\x$ denotes the lifting of the maps $h_\x: \ \T^* \to
\D$ (defined by (2.5)) to the ball $\Belt(f^{\mu_1}(D_\dl))_1$,  and
$\{P_n\}_1^\iy$ is a well-defined orthonormal polynomial basis in
$A_1^2(f^{\mu_1}(D_\dl))$ such that the degree of $P_n$ equals $n$
(canonically determined by the quasidisk $f^{\mu_1}(D_\dl)$; cf.
\cite{Mi}, \cite{Kr8}).

The quadratic differential $\psi_\dl$ in (3.8) has in the domain
$f^{\mu_1}(D_\dl)$ only zeros of even order, thus one can again
apply Proposition 2.1 with $D_0 = f^{\mu_1}(D_\dl)$ getting, similar
to (3.5), the equalities
  \be\label{3.9}
\tanh c_{\T^*}(\mathbf 0, S_{f^{\sigma_0}}) =
\vk_{f^{\mu_1}(D_\dl)}(f^{\mu_0}) = k = \tanh d_{\T^*}(\mathbf 0,
S_{f^{\sigma_0}}).
\end{equation}
On the other hand, since both \Ko \ and \Ca \ metric are
contractible under \hol \ maps and from (1.1),
$$
d_\T(\mathbf 0, S_{f^{\mu_0}}) = \tanh^{-1} k \ge c_\T(\mathbf 0,
S_{f^{\mu_0}}) \ge c_{\T^*}(\mathbf 0, \eta (S_{f^{\mu_0}})) =
c_{\T^*}(\mathbf 0, S_{f^{\sigma_0}}).
$$
Comparison with (3.9) implies
  \be\label{3.10}
c_\T(\mathbf 0, S_{f^{\mu_0}}) = d_\T(\mathbf 0, S_{f^{\mu_0}}) =
\tanh^{-1} k = \tau_\T(\mathbf 0, S_{f^{\mu_0}})
\end{equation}
proving the theorem in the case when one of the points is the origin
of $\T$.

The case of arbitrary two points $\vp = S_{f^\mu}, \ \psi =
S_{f^\nu}$ from $\T$ is investigated in a similar way (again by
applying Milin's coefficients), or can be reduced to (3.10) by the
right translations of type (3.6) moving one of these points to the
origin (a new base point of $\T$).

The proof for the infinitesimal metrics is similar. This completes
the proof of the theorem.

\section{Proof of Theorem 1.1}

First recall that the elements of the space $\T_1 = \T(\D^0)$ (where
$\D^0 = \D \setminus \{0\}$) are the equivalence classes of the \Be
\ coefficients $\mu \in \Belt(\D)_1$ so that the corresponding \qc \
automorphisms $w^\mu$ of the unit disk coincide on both boundary
components (unit circle $S^1 = \{|z| =1\}$ and the puncture $z = 0$)
and are homotopic on $\D \setminus \{0\}$. This space can be endowed
with a canonical complex structure of a complex Banach manifold and
embedded into $\T$ using uniformization.

Namely, the disk $\D^0$ is conformally equivalent to the factor
$\D/\G$, where $\G$ is a cyclic parabolic Fuchsian group acting
discontinuously on $\D$ and $\D^*$. The functions $\mu \in
L_\iy(\D)$ are lifted to $\D$ as the \Be \ $(-1, 1)$-measurable
forms  $\wt \mu d\ov{z}/dz$ in $\D$ with respect to $\G$, i.e., via
$(\wt \mu \circ \g) \ov{\g^\prime}/\g^\prime = \wt \mu, \ \g \in
\G$, forming the Banach space $L_\iy(\D, \G)$.

Extend these $\wt \mu$ by zero to $\D^*$ and consider the unit ball
$\Belt(\D, \G)$ of $L_\iy(\D, \G)$. Then the corresponding
Schwarzians $S_{w^{\wt \mu}|\D^*}$ belong to $\T$. Moreover, $\T_1$
is canonically isomorphic to the subspace $\T(\G) = \T \cap \B(\G)$,
where $\B(\G)$ consists of elements $\vp \in \B$ satisfying $(\vp
\circ \g) (\g^\prime)^2 = \vp$ in $\D^*$ for all $\g \in \G$. Most
of the results about the \uTs \ presented in Section 1 extend
straightforwardly to $\T_1$.

Due to the Bers isomorphism theorem, the space $\T_1$ is
biholomorphically equivalent to the Bers fiber space
$$
\mathcal F(\T) = \{\phi_\T(\mu), z) \in \T \times \C: \ \mu \in
\Belt(\D)_1, \ z \in w^\mu(\D)\}
$$
over the \uTs \ with \hol \ projection $\pi(\psi, z) = \psi$ (see
\cite{Be}). This fiber space is a bounded domain in $\B \times \C$.

\bk To prove the theorem, we establish the equalities (1.2) and
their infinitesimal counterpart for this fiber space.

We again model the space $\T$ as a domain in the space $\B$ formed
by the Schwarzians $S_{f^\mu}$ of functions $f^\mu(z) = z + b_0 +
b_1 z^{-1} + \dots \in \Sigma^0$ normalizing those additionally by
$f^\mu(1) = 1$.

Now the quadratic differentials defining the admissible \Te \
extremal coefficients $\mu_0 \in \Belt(\D)_1$ must be integrable and
\hol \ only on the punctured disk $\D \setminus \{0\}$ and can have
simple pole at $z = 0$, i.e., $\mu_0 = k |\psi_0|/\psi_0$ with
$$
\psi_0(z) = c_{-1} z^{-1} + c_0 + c_1 z + \dots \ , \quad 0 < |z| <
1.
$$
We associate with $f^\mu$ the odd function
 \be\label{4.1}
\mathcal R_{2,0} f^\mu(z) := (f^\mu(z^2) - f^\mu(0))^{1/2} = z +
\fc{ b_0 - f^\mu(0)}{2 z} + \fc {b_3^\prime}{z^3} + \dots
\end{equation}
whose \Gr \ coefficients $\a_{m n}(\mathcal R_{2,0} f^\mu)$ are
represented as polynomials of the first Taylor coefficients of the
original function $f^\mu$ and of $a = f^\mu(0)$. Hence, $\a_{m
n}(\mathcal R_{2,0} f^\mu)$ depend holomorphically on the
Schwarzians $\vp = S_{f^\mu} \in \T$ and on values $f^\mu(0)$, i.e.,
on pairs $X = (\vp, a)$ which are the points of the fiber space
$\mathcal F(\T)$. This joint holomorphy follows from Hartog's
theorem on separately \hol \ functions extended to Banach domains.
It allows us to construct for $\mathcal R_{2,0} f^\mu$ the
corresponding \hol \ functions (2.5) mapping the domain $\mathcal
F(\T)$ to the unit disk.

One can apply to $\mathcal R_{2,0} f^\mu$ the same arguments as in
the proof of Theorem 1.2 and straightforwardly establish for any \Te
\ extremal disk
$$
\{\phi_{\T_1}(t \mu_0) = X_t := (S_{f^{t\mu_0}}, f^{t\mu_0}(0)): \
|t| <1\} \quad (\mu_0 = |\psi_0|/\psi_0)
$$
in the space $\mathcal F(\T)$ the key equality
$$
\sup_{\x\in S(l^2)} \ |h_\x(S_{\mathcal R_{2,0} f^{t\mu_0}})| = |t|
$$
for any $|t| < 1$. This equality combined with (1.1) implies
 \be\label{4.2}
c_{\mathcal F(\T)}(\mathbf 0, X_t) = \tanh^{-1} |t| = \tau_{\mathcal
F(\T)}(\mathbf 0, X_t) = d_{\mathcal F(\T)}(\mathbf 0, X_t),
\end{equation}
and by Bers' biholomorphism between the spaces $\T_1$ and $\mathcal
F(\T)$, the similar equalities for the corresponding metrics on the
space $\T_1$. In view of the density of Strebel's points and
continuity of metrics, these equalities extend to all extremal disks
in $\T_1$, which yields the assertion of the theorem for the
distances of any point from the origin.

To establish the equality of distances between two arbitrary points
$X_1, \ X_2$ in $\T_1 = \T(D_{*})$, we uniformize the base point
$\D_{*} = \D \setminus \{0\}$ (with fixed homotopy class) by a
cyclic parabolic Fuchsian group $\G_0$ acting on the unit disk
(using the universal covering $\pi: \D \to \D_{*}$ with $\pi(0) =
0$) and embed the space $\T_1$ holomorphically into $\T$ via
$$
\T_1 = \T \cap \B(\D^*, \G_0) = \Belt(\D, \G_0)/\sim
$$
(where the equivalence relation commutate with the homotopy of \qc \
homeomorphisms of the surfaces). This preserve all invariant
distances on $\T_1$.

One can use the result of the previous step which provides that for
any point $X \in \T_1$ its distance from the base point $X_0 =
\D_{*}$ in any invariant (no-expanding) metric is equal to the \Te \
distance; hence,
$$
c_{\T_1} (X_0, X) = \tau_{\T_1} (X_0, X) = d_{\T_1} (X_0, X).
$$

Now, fix a \Be \ coefficient $\mu \in \Belt(\D, \G_0)_1$ so that
$X_1 = w^\nu(X_0)$ as marked surfaces (i.e., with prescribed
homotopy classes) and apply the change rule (3.6). It defines a \hol
\ automorphism $\sigma_\mu$ of the ball $\Belt(\D, \G_0)$ which is
an isometry in its \Te \ metric. This automorphism is compatible
with \hol \ factorizing projections $\phi_{\T_1}$ and
$\phi_{\T_1^*}$ defining the space $\T_1$ and its copy $\T_1^*$ with
the base point $X_1$. Thus $\sigma_\mu$ it descends to a \hol \
bijective map $\wh \sigma_\mu$ of the space $\T_1$ onto itself,
which implies the \Te \ isometry
$$
\tau_{\T_1} (\phi_{\T_1}(\mu), \phi_{\T_1}(\nu)) =
\tau_{\T_1}(\phi_{\T_1}(\mathbf 0), \phi_{\T_1}(\sigma_\mu(\nu)),
\quad \nu \in \Belt(f^\mu(\D), f^\mu \G_0 (f^\mu)^{-1})_1,
$$
and similar admissible isometries for the \Ca \ and \Ko \ distances
on these space. Combining this with the relations
$$
c_{\T_1}(\phi_{\T_1}(\mathbf 0), \phi_{\T_1}(\sigma_\mu(\nu))) =
\tau_{\T_1}(\phi_{\T_1}(\mathbf 0), \phi_{\T_1}(\sigma_\mu(\nu))) =
d_{\T_1}(\phi_{\T_1}(\mathbf 0), \phi_{\T_1}(\sigma_\mu(\nu)))
$$
established in the previous step, one derives the desired equalities
(1.2).

The case of infinitesimal metrics is investigated in a similar way,
which completes the proof of the theorem.

\section{Proofs of Theorem 1.3 and Corollary 1.5}

\subsection{Proof of Theorem 1.3} \ Theorems 1.1 and 1.2 imply,
together with the definition of complex geodesics, that these
geodesics in $\wt \T$ are the \Te \ geodesic disks
  \be\label{5.1}
\D(\mu_0) = \{\phi_{\wt \T}(t \mu_0/\|\mu_0\|_\iy): \ |t| < 1\}
\end{equation}
in this space. Accordingly, the uniqueness of the complex geodesic
joining a Strebel point in $\wt \T$ with the origin follows from
uniqueness of the disk (5.1) for such a point.

On the other hand, Tanigawa constructed in \cite{Ta} the extremal
\Be \ coefficients $\mu_0$ with nonconstant $|\mu_0(z)| <
\|\mu_0\|_\iy$ on a set of positive measure for which there exist
infinitely many distinct geodesic segments in the \uTs \ $\T$
joining the points $\phi_\T(\mathbf 0)$ and $\phi_\T(\mathbf
\mu_0)$. All these segments belong to different complex geodesics
joining the indicated points.

The case of geodesic disks joining two arbitrary points in $\wt \T$
is investigated in similar way. This completes the proof of the
theorem.

\bk\noindent {\bf Remark}. One can combine Theorems 1.1 and 1.2 with
the result of \cite{DTV} on existence of complex geodesics in convex
Banach domains and get an alternative proof of Theorem 1.3. The main
underlying facts ensuring the existence of geodesics are the
equality of invariant metrics established for geometrically convex
domains and weak$*$ compactness.

It is well known that if a Banach space $X$ has a predual $Y$, then
by the Alaougly-Bourbaki theorem the closure of its open unit ball
is weakly $^*$ compact. This holds, in particular, for our space
$\wt \T$ regarded as a bounded domain in $\B(\D^*, \G)$, which is
dual to the space $A_1(\D^*, \G)$ of integrable \hol \ quadratic
differentials with respect to group $\G$.

Theorems 1.1 and 1.2 ensure all the needed features, and therefore
one can obtain Theorem 1.3 also by applying the same arguments as in
\cite{DTV}.

\bk
\subsection{Proof or Corollary 1.5} We apply Golusin's improvement
of Schwarz's lemma which asserts that a \hol \ function
$$
g(t) = c_m t^m + c_{m+1} t^{m+1} + \dots: \D \to \D \quad (c_m \neq
0, \ \ m \ge 1),
$$
is estimated in $\D$ by
  \be\label{5.2}
|g(t)| \le |t|^m \fc{|t| + |c_m|}{1 + |c_m| |t|},
\end{equation}
and the equality occurs only for $g_0(t) = t^m (t+ c_m)/(1 + \ov c_m
t)$; see [Go, Ch. 8].

Fix $t_0 \ne 0$ and denote
$$
\mu_0^* = \mu_0/\|\mu_0\|_\iy, \quad \eta(t) = |t|^m (|t|+ |c_m|)/(1
+ |c_m| |t|).
$$
By Theorems 1.1 and 1.2, there exists a \hol \ map $j(\vp): \ \wt \T
\to \D$ (the limit \hol \ function for a maximizing sequence for the
\Ca \ distance) such that
 \be\label{5.3}
d_\D(0, |j \circ \phi_{\wt \T}(t_0 \mu_0^*)) = c_{\wt \T} (\mathbf
0, \phi_{\wt \T} (t_0 \mu_0^*)) = d_{\wt \T} (\mathbf 0, \phi_{\wt
\T} (t_0 \mu_0^*)).
\end{equation}
Thus the maps
$$
h(t) = \phi_{\wt \T}(t \mu_0^*): \D \to \D(\mu_0) \ \ \text{and} \ \
j_{*}(t) = j \circ \phi_{\wt \T}(t \mu_0^*) = j|\D(\mu_0): \
\D(\mu_0) \to \D
$$
determine two inverse hyperbolic isometries of the unit disk so that
$j_{*} \circ h(t) \equiv t$.

Now, let $J$ be a \hol \ functional on $\wt \T$ with the values in
$\D$ and its restriction $J_{\mu_0}$ to the disk $\D(\mu_0)$ is
expanded via (1.3). Then, using the relations (5.2) and (5.3) and
noting that $|\eta(t)| \le |t|$, one derives
$$
J_{\mu_0}(t_0) \le j_{*}(\eta(t_0)) = d_{\wt \T}(\mathbf 0,
h(\eta(t_0))) \le d_{\wt \T}(\mathbf 0, h(|t_0|))
$$
which implies (1.4). The case of equality easily follows from
Schwarz's lemma. This completes the proof of the corollary.

\section{Applications to geometric function theory}

\subsection{General distortion theorem}
The above theorems reveal the fundamental facts of the variational
theory for univalent functions with \qc \ extension.

Let again $L$ be a bounded oriented quasicircle in the complex plane
$\C$ separating the origin and the infinite point, with the interior
and exterior domains $D$ and $D^*$ so that $0 \in D$ and $\iy \in
D^*$. Put
$$
\Sigma_{k'}(D) = \{f \in \Sigma^0(D): \ k(f) \le k'\}.
$$

Consider on the class $\Sigma^0(D^*)$ a \hol \ (continuous and
Gateaux $\C$-differentiable) functional $J(f)$, which means that for
any $f \in \Sigma^0(D^*)$ and small $t \in \C$,
 \be\label{6.1}
J(f + t h) = J(f) + t J_f^\prime(h) + O(t^2), \quad t \to 0,
\end{equation}
in the topology of uniform convergence on compact sets in $\D^*$.
Here $J_f^\prime(h)$ is a $\C$-linear functional.

Assume that $J$ is lifted by $\wh J(\mu) = J(f^\mu)$ to a \hol \
function on $\Belt(D)_1$ and also depends holomorphically from the
Schwarzian derivatives $S_{f^\mu}$ on \uTs \ $\T$. Then the linear
functional $J_f^\prime(h)$ in (5.1) is the strong (Fr\'{e}chet)
derivative of $J$ in both norms $L_\iy$ and $\B(D^*)$.

Varying $f$ by (3.1), one gets the functional derivative
 \be\label{6.2}
\psi_0(z) = J_{\id}^\prime(g(\id, z))
\end{equation}
for the variational kernel
  \be\label{6.3}
g(w, \z) = 1/(w - \z) - 1/\z,
\end{equation}
Note that any such functional $J$ is represented by a complex Borel
measure on $\C$, which allows to extend this functional to all \hol
\ functions on $D^*$ (cf. \cite{Sc}). In particular, the value
$J_{\id}(g(\id, z))$ of $J$ on the identity map $\id(z) = z$ is
well-defined.

We assume that this derivative is meromorphic on $\C$ and has in the
domain $D$ only a finite number of the simple poles (hence $\psi_0$
is integrable over $D$). All this holds, for example, in the case of
the distortion functionals of the general form
$$
J(f) := J(f(a_1), \dots \ , f(a_m); \ f(z_1), f'(z_1), \dots \ ,
f^{(\a_1)}(z_1); \dots; f(z_p), f'(z_p), \dots \ , f^{(\a_p)}(z_p)).
$$
with $\wh J(\mathbf 0) = 0$ and $\grad \wh J(\mathbf 0) \ne 0$. Here
$a_1, \dots, \ a_m$ are distinct fixed points in $D$, and $z_1,
\dots \ , z_p$ are distinct fixed points in $D^*$ with assigned
orders $\a_1, \ \dots \ , \a_p$, respectively.

To have a possibility to apply Theorems 1.1 and 1.2, we restrict
ourselves by the model case $m = 1$, i.e., by the functionals
 \be\label{6.4}
J(f) = J(f(a); \ f(z_1), f'(z_1), \dots \ , f^{(\a_1)}(z_1); \dots;
f(z_p), f'(z_p), \dots \ , f^{(\a_p)}(z_p))
\end{equation}
depending on the values of maps at one point in the domain of
quasiconformality. In this case,
 \be\label{6.5}
\wh J_{\id}^\prime(g(\id, z)) = \fc{\partial \wh J(\mathbf
0)}{\partial \om} g(z, a) + \sum\limits_{j=1}^p
\sum\limits_{k=0}^{\a_j-1} \ \fc{\partial \wh J(\mathbf 0)}{\partial
\om_{j,k}} \fc{d^k}{d\z^k} g(w,\z)|_{w=z,\z=z_k},
\end{equation}
where $\om = f(a), \ \om_{j,k} = f^{(k)}(z_j)$; hence, $\psi_0$ is a
rational function.

For such functionals, Theorem 1.1 provides a general distortion
theorem which shed light on underlying features and, on the other
hand, implies the sharp explicit bounds.

\begin{thm} (i) \ For any functional $J$ of type (6.4) whose range domain
$J(\Sigma^0(D^*))$ has more than two boundary points, there exists a
number $\kp_0(J) > 0$ such that for all $\kp \le \kp_0(J)$, we have
the sharp bound
 \be\label{6.6}
\max\limits_{k(f) \le \kp} |J(f^\mu) - J(\id)| \le \max\limits_{|t|=
\kp} |J(f^{t|\psi_0|/\psi_0}) - J(\id)|;
\end{equation}
in other words, the values of $J$ on the ball $\Belt(D)_\kp = \{\mu
\in \Belt(D)_1: \ \|\mu\|_\iy \le \kp\}$ are placed in the closed
disk $\D(J(\id), M_\kp)$ with center at $J(\id)$ and radius $ M_\kp
= \max\limits_{|t|= \kp} |J(f^{t|\psi_0|/\psi_0}) - J(\id)|$. The
equality occurs only for $\mu = t|\psi_0|/\psi_0$ with $|t| = \kp$.

(ii) \ Conversely, if a functional $J$ is bounded via (6.6) for $0 <
\kp \le \kp_0(J)$ with some $\kp_0(J) > 0$, then up to rescaling
(multiplying $J$ by a positive constant factor),
 \be\label{6.7}
J(f^\mu) = g(S_{f^\mu}) + O(\|\mu\|_\iy^2) \quad \text{as} \ \
\|\mu\|_\iy \to 0,
\end{equation}
where $g$ is \hol \ on $\T_1$ and its renormalization $\wt g(\vp) =
g(\vp)/\sup_{\vp\in \T_1} |g(\vp)|$ is the defining map for the disk
$\D(\mu_0)$ as a $c_{\T_1}$-geodesic in the space $\T_1$ with the
base point representing the punctured quasidisk $D \setminus \{a\}$.
\end{thm}

\bk\noindent{\bf Outline of the proof}. First note that as one can
see from (6.7) that the underlying features arise from the
connection between such \hol \ functionals and the corresponding
$c_{\T_1}$-geodesics. The proof of Theorem 1.1 implies that the
restriction of $\wt g(\vp)$ to the disk formed by $\vp =
S_{f^{\mathcal R_{2,0}^* \mu_0}}$ is the inverse function for the
limit of the functions $h_\x(t \mathcal R_{2,0}^* \mu_0), \ |t| <
1$,  defined by (2.5).

The results of such type were obtained in \cite{Kr6}, \cite{Kr7} for
more specific functionals which relate to complex geodesics in the
\uTs \ $\T$. The proof of Theorem 6.1 involves $c_{\T_1}$-geodesics
and follows the same lines. Thus we only outline the main steps.

One can replace the assumption $f^\mu(0) = 0$ for $f^\mu \in
\Sigma^0(D)$ by $f^\mu(1) = 1$ and use the variation
 \be\label{6.8}
\begin{aligned}
f^\mu(z) &= z - \fc{1}{\pi} \iint\limits_D \mu(\z) \Bigl( \fc{1}{\z
- z} - \fc{1}{\z - 1} \Bigr) d \xi d \eta
+ O(\|\mu\|_\iy^2)  \\
&= z - \fc{z - 1}{\pi} \iint\limits_D \fc{\mu(\z) d \xi d \eta}{ (\z
- 1)(\z - z)} + O(\|\mu\|_\iy^2) \quad \text{as} \ \ \|\mu\|_\iy \to
0,
\end{aligned}
\end{equation}
hence replace (6.3) by
 \be\label{6.9}
g(w, \z) = \fc{1}{w - \z} - \fc{1}{w -1}.
\end{equation}
Assume also that $J(\id) = 0$, and let $f_0$ be a maximizing
function for $|J|$ in $\Sigma_\kp(D)$ (whose existence follows from
compactness). Take its extremal extension to $D$, i.e., with \Be \
coefficient
$$
\|\mu_{f_0}\|_\iy = \inf \{ \|\mu\|_\iy \le \kp: \ f^\mu = f_0 \ \
\text{on} \ \ D^* \cup \{a_1, \dots \ , a_m\}\},
$$
and suppose that
 \be\label{6.10}
\mu_{f_0} \ne \mu_0,
\end{equation}
where $\mu_0 = t |\psi_0|/\psi_0$ for some $t$ with $|t| = \kp$. Our
goal is to show that for small $\kp$ this leads to a contradiction.

Pick in $A_1(D \setminus \{a\})$  the functions
 \be\label{6.11}
\om_p(z) = \chi(z)^p \chi^\prime(z) - 1 - \psi_0(z), \quad p = 1, 2,
\dots \ ,
\end{equation}
where $\chi$ is a conformal map of $\D$ onto $D$ with $\chi(0) = 0,
\ \chi^\prime(0) > 0$ (hence, $\chi(z)^p \chi^\prime(z) = c^p z^p +
O(z^{p+1})$ as $z \to 0$), and
 \be\label{6.12}
\rho_a(z) = \fc{a - 1}{(z - 1)(z - a)}.
\end{equation}
For each $z \in D$, the function
$$
r_z(\z) = \fc{1}{\z - z} - \fc{1}{\z - 1}
$$
is in $A_1(D)$, hence $r_z(\z) = \sum\limits_0^\iy d_p(z) \z^p, \ \z
\in D$.

One of the main points in the proof is the following

\begin{lem} For sufficiently small $\kp \le \kp_0(J)$, the extremal \Be \
coefficient $\mu_{f_0}$ is orthogonal in $A_1(D)$ to all functions
(6.11) and (6.12), i.e., $\langle \mu_{f_0}, \rho_a\rangle_D = 0$
and $\langle \mu_{f_0}, \om_p\rangle_D = 0$ for all $p$.
\end{lem}

Its proof involves the properties of the projections of norm $1$ in
Banach spaces presented in \cite{EK} and investigation of norms
$$
h(\xi) = \iint\limits_\D |\psi_0(z) + \xi \rho_a(z)| dx dy, \quad
h_p(\xi) = \iint\limits_\D |\psi_0(z) + \xi \psi_p(z)| dx dy.
$$

Now apply a geodesic \hol \ map $g: \T_1 \to \D$ from Theorem 1.1
defining the disk $\D(\mu_0)$ as $c_{\T_1}$-geodesic; it determines
a hyperbolic isometry between this disk and $\D$. We lift this map
onto $\Belt(\D)_1$ by $\Ld(\mu) = g \circ \phi_{\T_1}(\mu)$ getting
a \hol \ map of this ball onto the disk. The differential of $\Ld$
at $\mu = \mathbf 0$ is a linear operator $P: L_\iy(\D) \to
L_\iy(\D)$ of norm $1$ which is represented in the form
$$
P(\mu) = \beta \langle \mu, \psi_0\rangle_D \ \mu_0.
$$
Let $P(\mu_{f_0}) = \a(\kp) \mu_0$. Since, by assumption, $f_0$ is
not equivalent to $f^{t_0\mu_0}$ with $|t_0| = \kp$, we have
$$
\Big\{\Ld \Bigl(\fc{t}{\kp}\mu_{f_0}\Bigr): \ |t| < 1\Big\}
\subsetneqq \{|t| < 1\}.
$$
Thus, by Schwarz's lemma,
 \be\label{6.13}
|\a(\kp)| < \kp.
\end{equation}
Note that the conjugate operator
$$
P^*(\psi) = \langle \mu_0, \psi\rangle_D \ \psi_0
$$
maps $L_1(D)$ into $L_1(D)$ and fixes the subspace $W = (\om_p,
\rho_a)$ of $A_1(D \setminus \{a\})$ spanned by functions (6.11) and
(6.12)

Now consider the function
 \be\label{6.14}
\nu_0 = \mu_{f_0} - \a(\kp) \mu_0
\end{equation}
which is not equivalent to zero, due to our assumption (6.10). Lemma
6.2 allows us to establish that $\nu_0$ annihilates all functions
from $\psi \in W$ and therefore orthogonal to all functions from the
whole space $A_1(D \setminus \{a\})$, because $\psi_0 \ \rho_a$ and
$\om_p, \ p = 1, 2, \dots$ form a complete set in this space. This
means that the function (6.14) belongs to the set
$$
A_1(D \setminus \{a\})^\bot = \{\mu \in L_\iy(D): \ \langle \mu,
\psi\rangle_D = 0 \ \ \text{for all} \ \ \psi \in A_1(D \setminus
\{a\})\}.
$$

But the well-known properties of extremal \qc \ maps imply that for
any $\nu \in A_1(D \setminus \{a\})^\bot$,
$$
\|\mu_{f_0}\|_\iy = \inf \{|\langle \mu_{f_0} + \nu, \psi \rangle_D:
\ \psi \in A_1(\D \setminus \{a\}), \ \|\psi\| = 1 \} \le
\|\mu_{f_0} + \nu\|_\iy.
$$
and therefore
$$
\|\mu_{f_0}\|_\iy = \kp \le \|\mu_{f_0} - \nu_0\|_\iy = \|\a(\kp)
\mu_0\|_\iy = \a(\kp),
$$
which contradicts (6.13). Hence $f_0$ is equivalent to
$f^{t|\psi_0|/\psi_0}$ and we can take $\mu_{f_0} = t
|\psi_0|/\psi_0$ for some $|t| = \kp$, completing the proof of the
first part of the theorem.

\bk To prove the converse assertion {\em (ii)}, we lift the original
functional $J$ to
$$
I(\mu) = \pi^{-1} \circ J(f^\mu): \ \Belt(D)_1 \to \D,
$$
where $\pi$ is a \hol \ universal covering of the domain $V(J) =
J(\Sigma^0(D^*))$ by a disk $\D_a = \{|z| < a\}$ with $\pi(0) = 0,
\pi^\prime(0) = 1$ (the lifting is single valued, since the ball
$\Belt(D)_1$ is simply connected). Let again $J(\id) = 0$. The
normalization of $\pi$ ensures that for sufficiently small $|\z|$,
$$
\pi(\z) = \z + O(\z^2)
$$
(with uniform estimate of the remainder for $|\z| < |\z_0|$), which
implies the asymptotic equality (6.7). The covering functional $I$
is \hol \ also in the Schwarzians $S_{f^\mu}$, which generates a
\hol \ map $\wt I: \ \T_1 \to \D$ so that $I = \wt I \circ
\phi_{\T_1}$. The above arguments provide for $I$ instead of (6.6)
the bound
$$
\max\limits_{k(f^\mu) \le \kp} |I(f^\mu)| = \kp \quad \text{for} \ \
0 < k < k_1(I).
$$
Restricting the covering map $\wt I$ to the extremal disk
$\{\phi_{\T_1}(t \mu_0^*): |t| < 1\} \subset \T_1$ (where $\mu_0^* =
|\psi_0|/\psi_0$) and applying to this restriction Schwarz's lemma,
one derives that $\wt I(\phi_{\T_1}(t \mu_0^*)) \equiv t$. Thus the
inverse to this map must be $c_\T$-geodesic, which completes the
proof of the theorem.

Representing the extremal $f^{t|\psi_0|/\psi_0}$ by (6.8), one can
rewrite the estimate (6.6) for $\kp \le \kp_1(J)$ in the form
  \be\label{6.15}
\max\limits_{k(f) \le \kp} |J(f^\mu) - J(\id)| \le \fc{\kp}{\pi}
\iint_D |J_{\id}^\prime(g(\id, z))| dx dy = \fc{\kp}{\pi}
\|\psi_0\|_1.
\end{equation}

\subsection{A lower estimate for the bound $\kp_0(J)$}. If the
functional $J$ is bounded on the whole class $\Sigma^0(D^*)$, and
$J(\id) = 0, \ \grad J(\id) \ne 0$, one can also derive from the
above arguments also a useful lower bound for $\kp_0(J)$, which
allows one to apply Theorem 6.1 effectively. Namely, one can verify
that the above proof works for
 \be\label{6.16}
\kp \le \kp_0(J) = \fc{\|J_{\id}^\prime\|}{\|J_{\id}^\prime\| + M(J)
+ 1},
\end{equation}
where
$$
\|J_{\id}^\prime\| = \fc{1}{\pi}\|\psi_0\|_1, \quad \quad M(J) =
\sup_{\Sigma^0(D^*)} |J(f)|.
$$

\subsection{ Additional remarks}

\bk\noindent{\bf 1.}  Similar theorem holds also for univalent
functions on bounded quasidisks $D$, for example, for the canonical
class $S_\kp(D)$ of univalent functions in $D$ normalized by $f(z) =
z + c_2 z^2 + \dots$ near the origin (provided that $z = 0 \in D$)
and admitting $\kp$-\qc \ extensions to $\hC$ which preserve the
infinite point. Such functions are investigated in the same way.

\bk\noindent{\bf 2.} Theorem 6.1 provides various explicit estimates
controlling the distortion in both conformal and \qc \ domains
simultaneously (comparing the known very special results established
in \cite{GR}, \cite{Kr1}, \cite{Ku1}, \cite{Ku2}).

\bk\noindent{\bf 3.}  The assumption that the distinguished point
$a$ is inner, is essential, and the estimate (6.6) can fail when the
functionals depend on values $f(a)$ at a prescribed point on the
boundary $\partial D$. The reasons are not technical. Actually the
bound $\kp_0(J)$ depends on the distance $\dist(a, \partial D)$ and
generically decreases to $0$ when $a$ approaches the boundary.

One can see this from the well-known result of K\"{u}hnau' on the
domain of values of $f(1)$ on $\Sigma_\kp$ presented, for example,
in [KK, Part 2]; it shows that in such a case an additional
remainder $O(\kp^2)$ can appear.

\section{New phenomena}

\subsection{Rigidity of extremals} The intrinsic connection between the
extremals of the distortion functionals on functions with \qc \
extensions and complex geodesics causes surprising phenomena which
do not appear in the classical theory concerning all univalent
functions. The differences arise from the fact that in problems for
the functions with \qc \ extensions the extremals belong to compact
subsets of $\Sigma^0(D^*)$ (or in other functional classes), while
the maximum on the whole class is attained on the boundary
functions.

We first mention the following consequence of Theorem 1.3 which
provides strong rigidity of extremal maps.

\begin{cor} In any class of univalent functions with $\kp$-\qc \
extension, neither function can be simultaneously extremal for
different \hol \ functionals (6.4) unless one of these functionals
have equal $1$-jets at the origin.
\end{cor}

\subsection{Example: the coefficient problem for functions with \qc \
extensions}
We mention here an improvement in estimating the
Taylor coefficients. Though the Bieberbach conjecture for the
canonical class $S$ of univalent functions  $f(z) = z +
\sum\limits_2^\iy a_n z^n$ in $\D$ has already been proved by de
Brange's theorem, the old coefficient problem remains open for
univalent functions in the disk with \qc \ extensions. The problem
was solved by the author for the functions with sufficiently small
dilatations.

Denote by $S_\kp(\iy)$ and $S_\kp(1)$ the classes of $f \in S$
admitting $\kp$-\qc \ extensions $\wh f$ to $\hC$ normalized by $\wh
f(\iy) = \iy$ and $\wh f(1) = 1$, respectively. Let
 \be\label{7.1}
 f_{1,t}(z) = \fc{z}{(1 - t z)^2}, \quad |z| < 1, \ \ |t| < 1.
\end{equation}
This function can be regarded as a \qc \ counterpart of the
well-known Koebe function which is extremal for many functionals on
$S$.

As a special case of Theorem 6.1, we have a complete solution of the
K\"{u}hnau-Niske problem \cite{KN} given by

\begin{thm} \cite{Kr3} For all $f \in S_\kp(\iy)$
and all $\kp \le 1/(n^2 + 1)$,
 \be\label{7.2} |a_n| \le 2 \kp/(n - 1),
\end{equation}
with equality only for the functions
  \be\label{7.3}
f_{n-1,t}(z) = f_{1,t}(z^{n-1})^{1/(n-1)} = z + \fc{2 t}{n - 1} z^n
+ \dots, \quad n = 3, 4, \dots; \ \ |t| = \kp.
\end{equation}
The estimate (7.2) also holds in the classes $S_k(1)$ with the same
bound for $\kp$.
\end{thm}

Note that every function (7.3) admits a \qc \ extension $\wh
f_{n-1,t}$ onto $\D^* = \{|z| > 1\}$ with \Be \ coefficient $t
\mu_n(z) = t |z|^{n+1}/z^{n+1}$, and $\wh f_{n-1,t}(\iy) = \iy$.

\bk No estimates have been obtained for arbitrary $\kp < 1$, unless
$n = 2$; in the last case, $|a_2| \le 2 \kp$ with equality for the
function (7.1) when $|t| = k$ (cf. \cite{Ku1}, \cite{KK},
\cite{KN}).

The rigidity provided by Corollary 7.1 yields that {\em the function
(7.1) cannot maximize $|a_n|$ in $S_\kp(\iy)$ even for one $\kp <
1$, unless} $n = 2$. Hence, for all $\kp < 1$,
 \be\label{7.4}
\max_{f\in S_\kp(\iy)} |a_n| > n \kp^{n-1} \quad (n \ge 3).
\end{equation}
For $n = 3$, this inequality was established in \cite{KN} involving
the elliptic integrals.

Comparing the coefficients $a_n$ of $f_{1,t}$ and $f_{n-1,t}$, one
derives from (7.2) and (7.4) the rough bounds for the maximal value
$\kp_n$ of admissible $\kp$ in (7.2):
$$
\fc{1}{n^2 + 1} \le \kp_n < \left[ \fc{2}{n(n -
1)}\right]^{1/(n-2)}.
$$

\subsection{Over-normalized functions}
Another remarkable thing in the distortion theory for univalent
functions with \qc \ extension concerns over-determined
normalization what reveals the intrinsic features of
quasiconformality. The variational problems for such classes are
originated in 1960s; the results were established mainly in terms of
inverse extremal functions $f_0^{-1}$ (see \cite{Kr1}, \cite{BK},
\cite{Re}).

We establish here some general explicit bounds. Assume that $z = 1$
lies on the common boundary of $D$ and $D^*$ which separates the
points $0$ and $\iy$ and denote by $\Sigma^0(D^*, 1)$ the class of
univalent functions in $D^*$ with \qc \ extensions across $L$ which
satisfy
$$
f(z) = z + \const + O(1/z) \ \ \text{near} \ \ z = \iy; \ \ f(1) =
1,
$$
and by $\Sigma_\kp(D^*, 1)$ its subclasses consisting of functions
with $\kp$-\qc \ extensions. Fix in the complementary domain $D$ a
finite collection of points
$$
e = (e_1, \dots \ , e_m)
$$
and associate with this set the following subspaces of $L_1(D)$: the
span $\mathcal L(e)$ of rational functions
$$
\rho_s(z) =  \fc{e_s - 1}{(z - 1)(z - e_s)} \quad s = 1, \dots , m,
$$
the space $A_1(D_e)$ of integrable \hol \ functions in the punctured
domain $D_e = D \setminus \{e_1, \dots \ , e_m\}$, and
$$
\mathcal L_0 = \mathcal L(e) \bigoplus \{c \psi_0: \ c \in \C\},
$$
where $\psi_0 = J_{\id}^\prime(g(\id, \cdot))$ for $g(w, \z)$ given
by (6.8). Let
  \be\label{7.5}
\Sigma_\kp(D^*, 1, e) = \{f \in \Sigma_\kp(D^*, 1): \ f(e_s) = e_s,
\ \ s = 1, \dots , m\}, \quad \Sigma^0(D^*, 1, e) = \bigcup_{\kp}
\Sigma_\kp(D^*, 1, e).
\end{equation}

Note that these classes with over-determined normalization contain
nontrivial maps $f^\mu \neq \id$ for any $\kp < 1$ what is insured,
by the local existence theorem from [Kr11, Ch. 4]. We shall use its
special case for simply connected plain domains presenting it as

\begin{lem}Let $D$ be a simply connected domain on the Riemann sphere
$\hC$. Assume that there are a set $E$ of positive two-dimensional
Lebesgue measure and a finite number of points
 $z_1, z_2, ..., z_m$ distinguished in $D$. Let
$\a_1, \a_2, ..., \a_m$ be non-negative integers assigned to $z_1,
z_2, ..., z_m$, respectively, so that $\a_j = 0$ if $z_j \in E$.

Then, for a sufficiently small $\ve_0 > 0$ and $\varepsilon \in (0,
\varepsilon_0)$, and for any given collection of numbers $w_{sj}, s
= 0, 1, ..., \a_j, \ j = 1,2, ..., m$ which satisfy the conditions
$w_{0j} \in D$, \
$$
|w_{0j} - z_j| \le \ve, \ \ |w_{1j} - 1| \le \ve, \ \ |w_{sj}| \le
\ve \ (s = 0, 1, \dots   a_j, \ j = 1, ..., m),
$$
there exists a \qc \ self-map $h$ of $D$ which is conformal on $D
\setminus E$ and satisfies
$$
h^{(s)}(z_j) = w_{sj} \quad \text{for all} \ s =0, 1, ..., \a_j, \ j
= 1, ..., m.
$$
Moreover, the \Be \ coefficient $\mu_h(z) = \partial_{\bar z}
h/\partial_z h$ of $h$ on $E$ satisfies $\| \mu_h \|_\iy \leq M
\ve$. The constants $\ve_0$ and $M$ depend only upon the sets $D, E$
and the vectors $(z_1, ..., z_m)$ and $(\a_1, ..., \a_m)$.

If the boundary $\partial D$ is Jordan or is $C^{l + \a}$-smooth,
where $0 < \a < 1$ and $l \geq 1$, we can also take $z_j \in
\partial D$ with $\a_j = 0$ or $\a_j \leq l$, respectively.
\end{lem}

Let us estimate on such over-normalized classes the functionals
 \be\label{7.6}
J(f) = J(f(z_1), f'(z_1), \dots \ , f^{(\a_1)}(z_1); \dots; f(z_p),
f'(z_p), \dots \ , f^{(\a_p)}(z_p))
\end{equation}
controlling the distortion on the domain of conformality.

Now one can use only conditional \qc \ variations whose \Be \
coefficients are orthogonal to the rational quadratic differentials
(6.12) corresponding to the fixed points. Thus the above proof of
key Lemma 6.2 fails, and this Lemma and Theorem 6.1 do not work.

The following theorem provides the sharp explicit bounds for
sufficiently small $\kp$ involving $L_1$-distance between the
functional derivative $\psi_0$ and span $\mathcal L(e)$.

\begin{thm} For any functional (7.5) and any finite set $e$ of fixed
points in $D$, there exists a positive number $\kp_0(J, e) < 1$ such
that for all $\kp \le \kp_0(J, e)$, we have for any function $f \in
\Sigma_\kp(D^*, 1, e)$  the sharp bound
 \be\label{7.7}
\max\limits_{\|\mu\|\le \kp} \ |J(f^\mu) - J(\id)| = |J(f^{\kp
|\psi_e|/\psi_e}) - J(\id)| = d \kp + O(\kp^2)
\end{equation}
with uniformly bounded ratio $O(\kp^2)/\kp^2$, where
 \be\label{7.8}
\psi_e = \xi_0 \psi_0 + \sum\limits_1^m \xi_s \rho_s
\end{equation}
with some constants $\xi_0, \xi_1, \dots, \ \xi_m$ and
 \be\label{7.9}
d = \inf_{\mathcal L(e)} \ \|\xi_0\psi_0 - \psi\|_1.
\end{equation}
The constants $\xi_0, \xi_s$ in (7.8) are determined (not necessary
uniquely) by the conditions
 \be\label{7.10}
\langle |\psi_e|/\psi_e, \psi \rangle_D = 0 \quad \text{for all} \ \
\psi \in \mathcal L(e); \ \ \langle |\psi_e|/\psi_, \psi_0 \rangle_D
= d.
\end{equation}
\end{thm}

\bk\noindent{\bf Proof}. By the Hahn-Banach theorem, there exists a
linear functional $l$ on $L_1(D)$ such that
  \be\label{7.11}
l(\psi) = 0, \ \ \psi \in \mathcal L(e); \ \ l(\psi_0) = d,
\end{equation}
and
$$
\|l\|_{L_1(D)} = \|l\|_{\mathcal L_0} = 1,
$$
and this norm is minimal on the spaces $\mathcal L_0 \subset A_1(D
\setminus e) \subset L_1(D)$. Hence, for any other linear functional
$\wt l$ on $L_1(D)$ satisfying (7.9) must be $ \|\wt l\|_{\mathcal
L_0} \ge 1$ and $\wt l(\psi_0) \le d$; otherwise, were $\wt
l(\psi_0) = r d$ with $r > 1$, the functional $\wt l/r$ with norm
less than $1$ would satisfy (7.11), in contradiction to minimality.

The functional $l$ is represented on $L_1(D)$ via
$$
l(\psi) = \iint\limits_D \nu_0(z) \psi(z) dx dy , \quad \psi \in
L_1,
$$
with some $\nu_0 \in L_\iy(D)$ so that
  \be\label{7.12}
\iint\limits_D \nu_0(z) \psi(z) dx dy = 0, \ \ \psi \in \mathcal
L(e); \ \ \iint\limits_D \nu_0(z) \psi_0(z) dx dy = d.
\end{equation}
Since the norm of $l$ on the widest space $L_1(D)$ is attained on
its subspace $\mathcal L_0$, the function $\nu_0$ is of the form
$\nu_0(z) = |\psi_e(z)|/\psi_e(z)$ with integrable \hol \ $\psi_e$
on $D \setminus e$ given by (7.8).

After extending $\nu_0$ by zero to $D^*$, which yields  an extremal
\Be \ coefficient $\nu_0 \in \Belt(D)_1$ for our \hol \ functional
$J$, one can represent the map $f^{t\nu_0}$ by (6.8) getting from
the second equality in (7.12) and from indicated minimality of
$\|l\|$ the estimate (7.10), since for all other $f^\nu \in
\Sigma_{\kp}(D^*, 1, e)$, we have $|J(f^\nu)| \le |J(f^{\kp
\nu_0})|$.

However, this $\kp$-\qc \ map can move the fixed points $e_s$ to
$f_0(e_s) = e_s + O(\kp^2)$, where $f_0 = f^{\nu_0}$. Thus one needs
to apply additional $O(\kp^2)$-\qc \ variation $h_0$ by Lemma 7.3 to
get $h_0 \circ f_0(e_s) = e_s$ (for all $s$) and preserving the
values $f_0^{(\a_j)}(z_j)$, and then take the extremal map $\wt f$
(with smallest dilatation) satisfying
 \be\label{7.13}
\wt f^{(\a_j)}(z_j) = f_0^{(\a_j)}(z_j), \quad \wt f(e_s) = e_s
\end{equation}
for all given $\a_j$ and $e_s$ so that its defining \hol \ quadratic
differential $\psi_e$ belongs to the subspace $\mathcal L_0$. It can
be shown, using the uniqueness of \Te \ extremal maps generated by
integrable \hol \ quadratic differential that this $psi_e$ is unique
in $A_1(D \setminus e)$.

The assertion on uniform bound for the remainder in (7.10) follows
the general distortion results for \qc \ maps. This completes the
proof of the theorem.

\bk\noindent{\bf Remarks}.

{\bf 1}. The quadratic differential $\psi_e$ constructed in the
proof depends also on $\kp$.

{\bf 2}. The assumption $f(1) = 1$ can be replaced by $f(0) = 0$;
then the fixed points $e_s$ must be chosen to be distinct from the
origin.

\bk Similar theorem also holds for the over-normalized functions in
bounded quasidisks. We illustrate it on the coefficient problem:

{\em Find $\max |a_n| \ (n \ge 2)$ for the functions $f \in
S_\kp(\iy)$ leaving a given set $e = (e_1, \dots, \ e_m) \subset \D
\setminus \{0\}$ fixed}.

In this case, one gets as a consequence the following result.

\bk
\begin{thm}
For any $n \ge 2$, there is a number $\kp_n(e) < 1$ such that for
$kp \le \kp_n$ and all $f \in S_\kp(\iy)$, which fix a given set $e
= (e_1, \dots, \ e_m)$, we have the sharp bound
 \be\label{7.14}
\max\limits_{\|\mu\|\le \kp}  |a_n(f^\mu)| =
|a_n(f^{\kp|\psi_n|/\psi_n})| = d_n \kp + O(\kp^2),
\end{equation}
where similar to (7.8) and (7.9),
$$
\psi_n(z) = c z^{-n-1} + \sum\limits_1^m \xi_s \rho_s(z), \quad d_n
= \inf_{\mathcal L(e)} \|\psi_n - \psi\|_1.
$$
The constants $c, \xi_s$ are determined from the equations of type
(7.10), and the remainder in (7.14) is estimated uniformly for all
$\kp \le \kp_n$.
\end{thm}

\bk The distortion bounds of type (7.7) given by Theorem 7.4 and its
corollaries hold in somewhat weakened form (up to terms $O(\kp^2)$
for the maps preserving an infinite subset $e$ in $D$, provided that
the corresponding class $\Sigma_\kp(D^*, 1, e)$ contains the
functions $f^\mu \neq \id$.

The proof is similar but now the quadratic \hol \ differentials
$\psi_e$ defining the extremal functions are represented instead of
(7.8) in the form
$$
\psi_e = c \psi_0 + \psi, \quad \psi \in \mathcal L(e)
$$
and there are no variations of type Lemma 7.3 for the infinite sets.

\bk
\medskip
{\small\em{ \leftline{Department of Mathematics, Bar-Ilan
University} \leftline{5290002 Ramat-Gan, Israel} \leftline{and
Department of Mathematics, University of Virginia,}
\leftline{Charlottesville, VA 22904-4137, USA}}}


\bigskip
\bigskip
\begin{thebibliography}{EKK}
{\small

\bibitem[Ah] {Ah}
L.V. Ahlfors, {\em On \qc \ mappings},
J. Anal. Math. \textbf{3} (1953-1954), 1-58.

\bibitem[Be] {Be}
L. Bers, {\em Fiber space over \Te \ spaces},
Acta Math. \textbf{130} (1973), 89-126.

\bibitem[BK] {BK}
P. Biluta and S.L. Krushkal, {\em On the question of extremal
quasiconformal  mappings}, Soviet Math. Dokl. \textbf{11} (1971),
76-79.

\bibitem[Di] {Di}
S. Dineen, {\em The Schwarz Lemma},
Clarendon Press, Oxford, 1989.

\bibitem[DTV] {DTV}
S. Dineen, R.M. Timoney and J.P. Vigu\'{e},
{\em Pseudodistances invariantes sur les domains d'une espace
localement convexe},
Ann. Scuola. Norm. Sup. Pisa Cl. Sci.(4) \textbf{12} (1985), 515-529.

\bibitem[EE] {EE}
C.J. Earle and J.J. Eells, {\em On the differential geometry of \Te
\ spaces}, J. Analyse Math. \textbf{19} (1967), 35-52.

\bibitem[EK] {EK}
C.J. Earle and I. Kra, {\em On sections of some \hol \ families of
closed Riemann surfaces}, Acta Math. \textbf{137} (1976), 49-79.

\bibitem[EKK] {EKK}
C.J. Earle, I. Kra and S. L. Krushkal,
{\em Holomorphic motions and \Te \ spaces},
Trans. Amer. Math. Soc.  \textbf{343} (1994),   927-948.

\bibitem[EM] {EM}
C.J. Earle and S. Mitra,
{\em Variation of moduli under \hol \ motions},
In the tradition of Ahlfors and Bers (Stony Brook, NY, 1998),
Contemp. Math. \textbf{256}, Amer. Math. Soc., Providence, RI, 2000,
pp. 39-67.

\bibitem [GL] {GL}
F.P. Gardiner and N. Lakic, {\em Quasiconformal \Te \ Theory}, Amer.
Math. Soc., 2000.

\bibitem[Go] {Go}
G.M. Goluzin, {\em Geometric Theory of Functions of Complex
Variables}, Transl. of Math. Monographs, vol. 26, Amer. Math. Soc.,
Providence, RI, 1969.

\bibitem[Gr] {Gr}
H. Grunsky, {\em Koeffizientenbedingungen f\"{u}r schlicht abbildende
meromorphe Funktionen}, Math. Z. \textbf{45} (1939), 29-61.

\bibitem[GR] {GR}
V.Ya. Gutlyanskii and V.I. Ryazanov, {\em Geometric and Topological
Theory of Functions and Mappings (Geometricheskaya i
topologicheskaya teoria funkcii i otobrazhenii)}, Naukova Dumka,
Kiev, 2011 (Russian).

\bibitem[K] {K}
I. Kra, {\em The \Ca \ metric on abelian \Te \ disks}, J. Anal.
Math. \textbf{40} (1981), 129-143.

\bibitem[Kr1] {Kr1}
S.L. Krushkal, {\em Quasiconformal Mappings and Riemann Surfaces},
Wiley, New York, 1979.

\bibitem[Kr2] {Kr2}
S.L. Krushkal,
{\em Grunsky coefficient inequalities, Carath\'{e}odory metric and
extremal \qc \ mappings},
Comment. Math. Helv. \textbf{64} (1989), 650-660.

\bibitem[Kr3] {Kr3}
S.L. Krushkal, {\em Exact coefficient estimates for univalent functions
with quasiconformal extension},
Ann. Acad. Sci. Fenn. Ser. A.I. Math. \textbf{20} (1995), 349-357.

\bibitem[Kr4] {Kr4}
S.L. Krushkal {\em Complex geometry of the universal Teichm\"{u}ller
space},
Siberian Math. J. {\bf 45}(4) (2004), 646-668.

\bibitem[Kr5] {Kr5}
S.L. Krushkal  {\em Plurisubharmonic features of the \Te \
metric}, Publications de l'Institut Math\'{e}matique-Beograd,
Nouvelle s\'{e}rie \textbf{75(89)} (2004), 119-138.

\bibitem[Kr6] {Kr6}
S.L. Krushkal, {\em Extremal problems for Fredholm eigenvalues},
Israel J. Math. \textbf{172} (2009), 279-307.

\bibitem[Kr7] {Kr7}
S.L. Krushkal, {\em Hyperbolic metrics on \uTs \ and extremal problems},
Ukrains'kii Matematychnii Visnyk \textbf{8} (2011), 557-579;
J. of Math. Sciences \textbf{182}, no. 1 (2012), 70-86.

\bibitem[Kr8] {Kr8}
134. S.L. Krushkal, {\em Generalized Grunsky coefficient
inequalities and quasiconformal deformations}, Uzbek Math. J., 2014,
no. \textbf{1}, 30-45.

\bibitem[KK] {KK}
S. L. Kruschkal und R. K\"{u}hnau,
{\em  Quasikonforme Abbildungen - neue Methoden und Anwendungen},
 Teubner-Texte zur Math., vol. \textbf{54},  Leipzig, 1983.

\bibitem[Ku1] {Ku1}
R. K\"{u}hnau, {\em Verzerrungss\"atze und Koeffizientenbedingungen
vom Grunskyschen Typ f\"{u}r quasikonforme Abbildungen}, Math. Nachr.
\textbf{48} (1971), 77-105.

\bibitem[Ku2] {Ku2}
R. K\"{u}hnau, {\em \"{U}ber die Werte des Doppelverh\"{a}ltnisses
bei quasikonformer Abbildung},
Math. Nachr. \textbf{95} (1980), 237-251.

\bibitem[Ku3] {Ku3}
R. K\"{u}hnau,
{\em Wann sind die Grunskyschen Koeffizientenbedingungen hinreichend
f\"{u}r $Q$-quasikonfor\-me Fortsetzbarkeit}?
Comment. Math. Helv. \textbf{61} (1986), 290-307.

\bibitem[KN] {KN}
   R. K\"{u}hnau und W. Niske,
{\em Absch\"{a}tzung des dritten Koeffizienten bei den quasikonform
fortsetzbaren schlichten Funktionen  der Klasse $S$}, Math. Nachr.
\textbf{78} (1977), 185-192.

\bibitem[Le] {Le}
N.A. Lebedev, {\em The Area Principle in the Theory of Univalent
Functions}, Nauka, Moscow, 1975 (Russian).

\bibitem[Mi] {Mi}
I.M. Milin, {\em Univalent Functions and Orthonormal Systems},
Transl. of mathematical monographs, vol. 49,
Transl. of Odnolistnye funktcii i normirovannie systemy,
Amer. Math. Soc., Providence, RI, 1977.

\bibitem[Po] {Po}
Chr. Pommerenke, {\em Univalent Functions}, Vandenhoeck $\&$
Ruprecht, G\"{o}ttingen, 1975.

\bibitem[Re] {Re}
H. Renelt, {\em Extremalprobleme bei quasikonformen Abbildungen
unter h\"{o}heren Normierungen}, Math. Nachr. \textbf{66} (1975),
125-143.

\bibitem[Ro] {Ro}
H.L. Royden, {\em Automorphisms and isometries of \Te \ space},
Advances in the Theory of Riemann Surfaces (Ann. of Math. Stud.,
vol. 66), Princeton Univ. Press, Princeton, 1971, pp. 369-383.

\bibitem[SS] {SS}
M. Schiffer and D. Spencer, {\em Functionals of finite Riemann
Surfaces}, Princeton Univ. Press, Princeton, 1954.

\bibitem[Sc] {Sc}
G. Schober, {\em Univalent Functions - Selected Topics},
Lecture Notes in Math. \textbf{478}, Springer, Berlin, 1975.

\bibitem[St] {St}
K. Strebel, {\em On the existence of extremal Teichmueller
mappings}, J. Anal. Math. \textbf{30} (1976), 464-480.

\bibitem[Ta] {Ta}
H. Tanigawa, {\em Holomorphic families of geodesic discs in infinite
dimensional Teichm\"uller spaces}, Nagoya Math. J. \textbf{127}
(1992), 117-128.

\bibitem[Ve] {Ve}
E. Vesentini, {\em Complex geodesics and \hol \ mappings},
Sympos. Math. \textbf{26} (1982), 211-230.

\bibitem[Zh] {Zh}
I.V. Zhuravlev, {\em Univalent functions and \Te \ spaces},
Inst. of Mathematics, Novosibirsk, preprint, 1979, 1-23 (Russian).



}
\end{thebibliography}
\end{document}